\documentclass[letterpaper,oneside,english]{amsart}
\usepackage{amsmath}
\usepackage{amsthm} 
\usepackage{mathptmx}

  \newtheorem{lem}{Lemma}
   \newtheorem{thm}{Theorem}

\usepackage{latexsym,amssymb}

\newcommand{\U}{\Upsilon_{\omega}}
\newcommand{\E}{\mathcal{E}}
\newcommand{\R}{\mathbb{R}}
\newcommand{\Z}{\mathbb{Z}}
\newcommand{\Dom}{\operatorname{Dom}}
\newcommand{\supp}{\operatorname{supp}}
\newcommand{\Tr}{\operatorname{Tr}}
\newcommand{\del}{\partial}
\newcommand{\bdry}{\partial}
\newcommand{\grad}{\nabla}

\begin{document}
\title{The rate of decay of the Wiener sausage in local Dirichlet space.}
\author{Lee R. Gibson \and Melanie Pivarski}
\address{Lee R. Gibson, Department of Mathematics,
University of Louisville, Louisville, KY 40292 \\
Tel.: +1-502-852-6826, Fax: +1-502-852-7132 \\ 
lee.gibson@louisville.edu
and
Melanie Pivarski,  
Department of Mathematics and Actuarial Science,
Roosevelt University, 
430 S. Michigan Ave,
Chicago, IL 60605\\
mpivarski@roosevelt.edu
}
\date{28 July, 2010}
\begin{abstract}
In the context of a heat kernel diffusion which admits a Gaussian type estimate with parameter $\beta$ on a local Dirichlet space, we consider the log asymptotic behavior of the negative exponential moments of the Wiener sausage.  We show that the log asymptotic behavior up to time $t^{\beta}V(x,t)$ is $V(x,t)$, which is analogous to the Euclidean  result.  Here $V(x,t)$ represents the mass of the ball of radius $t$ about a point $x$ of the local Dirichlet space.  The proof uses a known coarse graining technique to obtain the upper asymptotic, but must be adapted for use without translation invariance in this setting.  This result provides the first such asymptotics for several other contexts, including diffusions on complete Riemannian manifolds with non-negative Ricci curvature.
\end{abstract}

\maketitle

\section{INTRODUCTION} 
\label{intro}
In the celebrated 1975 paper by Donsker and Varadhan, the negative exponential moment with parameter $\nu$ of the Wiener sausage in $\mathbb{R}^d$ was shown to be log asymptotically equivalent to $ -k(\nu,d)t^{d/(d+2)}$, for a known value of $k(\nu,d)$. \cite{DV}  Later, Sznitman's method of enlargement of obstacles was used to recover this result and to solve many other problems from the field of random media. \cite{SznBook} More recent work in the area of analysis on metric spaces has seen the development of local Dirichlet spaces as a very general context in which diffusion processes make sense in a natural way. \cite{LSCBendikov},\cite{LSCHebisch}, \cite{Sturm}, As knowledge of the properties of these spaces grows, it is natural to investigate the extent to which the Donsker and Varadhan result might hold for the diffusions that inhabit them.  As a first step into this investigation, we consider a local Dirichlet space $X$ satisfying some volume estimates on which the heat kernel admits a two-sided Gaussian type estimate with parameter $\beta \le 2$.  We show that for the sample path $C_t$ of the diffusion process (related to the heat kernel \cite{SturmDiffHK}) starting from $x \in X,$  
$$  \log \mathbb{E}^x\left[\exp\left(-\nu \mu\left(C_{t^{\beta}V(x,t)}^\epsilon\right)\right) \right] \simeq V(x,t),
$$
where $\nu>0$, $V(x,t)$  represents the volume of the ball of radius $t$ centered at $x\in X$, and $C_t^{\epsilon}$ is the union of the radius $\epsilon$ balls centered at the points of the sample path up to time $t$.  Here $f(t)\simeq g(t)$ means that there are positive constants $c$ and $C$ for which $cg(t)\le f(t) \le Cg(t)$ for all sufficiently large $t$.
This result may be applied directly to the examples of diffusions on the fractal-like manifolds and graphs which are described in \cite{BBK}, on complete Riemannian manifolds with non-negative Ricci curvature \cite{LiYau}, and to those on Riemannian complexes with an underlying finitely generated polynomial growth group \cite{PivarskiLSC}.

Since the proofs related to the known behavior of the Wiener sausage in $\mathbb{R}^d$ rely on the translation and scale invariance of $\mathbb{R}^d$ and Brownian motion, the techniques of those proofs cannot be applied directly in the local Dirichlet space context.  Instead, we will follow the solution technique of the discrete analogue of our problem (random walk on a non-transitive graph) appearing in \cite{Gibson}.  In this process, we will make use of the relationship between the exponential moments of the volume of the sausage and the annealed survival probability of the diffusion process among Poisson distributed hard obstacles.   In particular, the main result follows from the connection between this survival time and the probability of an abnormally large obstacle-free region.

\section{Assumptions}\label{Assumptions}
We would like to look at a space which has a sufficiently nice structure to allow for a Dirichlet form leading to a heat kernel.  We will require some bounds on both the heat kernel and the local and global geometry of the space.  

\subsection{Differentiable structure on the space} \label{Differentiable} We begin by considering a complete space $X$ which is a local Dirichlet space as described in \cite{Sturm}.  This is a metric measure space, $(X,\mu,d)$, with
a fixed regular, strongly local, symmetric Dirichlet form, $\E$, whose
domain, $\Dom(\E)$, is on the real Hilbert space $L^2(X,\mu)$ with
norm $||f||_2 = \sqrt{\int_X f^2 d\mu}$. Because we will later define a Poisson point process on $X$, the positive radon measure $\mu$ will be assumed to be $\sigma$-finite.  The Dirichlet form, $\E$, has an associated nonnegative semi-definite self-adjoint operator
on $L^2(X, \mu)$.  This operator is a natural analogue of the Laplacian; to emphasize this, we will write it as $\Delta$.   The following holds for all
$f\in \Dom(\Delta), g \in \Dom(\E)$:
\begin{eqnarray*}
\E(f,g) = \langle \Delta f,g \rangle = \int_X \Delta f(x)g(x)d\mu(x)
\end{eqnarray*}Because $\Delta$ is self-adjoint, there exists an operator
$\sqrt{\Delta}$ which acts on $f,g \in \Dom(\Delta)$ by:
\begin{eqnarray*}
\langle \Delta f,g\rangle = \langle \sqrt{\Delta} f, \sqrt{\Delta}g \rangle.
\end{eqnarray*}
This is an analogue of the gradient; we will use the notation $\grad =
\sqrt{\Delta}$.
The structure of a local Dirichlet space allows for a well-defined heat kernel \cite{SturmDiffHK}, \cite{SturmDif}.   The heat kernel, $h_t(x,y)$, is both the kernel of the semigroup
$H_t=e^{-t\Delta}$ and the fundamental solution to the heat
equation  $ \del_t u = - \Delta u.$  Our arguments will primarily rely upon the semigroup interpretation.  Note that we will not go about the often messy process of constructing such a kernel ourselves; rather, we will consider a space where the kernel exists and satisfies certain bounds.

At times, it will be useful to consider the operator $\Delta$ restricted to a connected open set $U \subset X$.  We denote this operator by $\Delta^U$.  In this context, the
smallest non-zero eigenvalue of $\Delta^{U}$ satisfies
\begin{align}
\lambda(U) & 
=\inf\left\{ \left.
\frac{\left\langle \Delta^{U}f,f\right \rangle}{\left\langle f,f\right\rangle}
\right|f\not\equiv0,\mathrm{supp}f\subset U\right\}. \nonumber 
\label{eq:eigenvalue}\end{align}
For such an operator, there is a corresponding semigroup $H_t^U$ and killed diffusion with kernel $h^U_t(x,y)$, which can be computed in in terms of the ordinary diffusion using the Dynkin-Hunt formula:
\[
h_t^{U}(x,y)=h_t(x,y) -\mathbb{E}^x(h_{t-\tau_U}(X_{\tau_U},y)1_{\tau_U\le t}).
\]

\subsection{Volume estimates} \label{Volume} We require our space to have comparable volumes at all points in the space $X$ and to grow in volume at least linearly with respect to the radius.  For this we consider balls of radius $r$ centered at $x$, denoted $B(x,r)$ with volume denoted $V(x,r) = \mu(B(x,r))$.  We require for some $c>0$ and all $r>0$
\begin{eqnarray}
cr < \inf_x V(x,r) \le \sup_x V(x,r) < \infty .\label{gb1}
\end{eqnarray}
Additionally, we require the following bound on the relative sizes of nearby balls as their radii increase:
\begin{equation}
\liminf_{r\to\infty}V(x,r)^{-1}\inf_{y\in B\left(x,r^{\beta}V(x,r)\right)}V(y,r)>0 .\label{eq:volcondbeta}\end{equation}
Note that such an estimate will hold whenever $V(x,r)\simeq f(r)$.

\subsection{Gaussian and sub-Gaussian heat kernel bounds} \label{Gaussian} We consider only spaces where the heat kernel satisfies a Gaussian-style estimate with parameter $\beta$.  We denote this estimate by GE($\beta$).  Such an estimate is satisfied exactly when there exist constants $c, C$ such that for all $x,y\in X$ and $t>0$ we have both the upper estimate \begin{equation}
  h_{t}(x,y)\le\frac{C}{V\left(x,t^{1/\beta}\right)}\exp\left(-c\left(\frac{d(x,y)^{\beta}}{t}\right)^{1/(\beta-1)}\right),\label{eq:GUE}\end{equation}
and the lower estimate\begin{equation}
  h_t(x,y)\ge\frac{c}{V\left(x,t^{1/\beta}\right)}\exp\left(-C\left(\frac{d(x,y)^{\beta}}{t}\right)^{1/(\beta-1)}\right).\label{eq:GLE}\end{equation}
Note that when $\beta=2$, this is the standard Gaussian bound that is satisfied by the heat kernel in $\R^d$.  
\subsection{Consequences of the Gaussian estimates}\label{Consequences}
The GE($\beta$) estimate implies both volume doubling and an associated  Poincar\'{e} inequality \cite{SturmAnalysisOn}, \cite{LSCHebisch}.  We describe both below.  

The space $X$ satisfies the Neumann Poincar\'{e}
inequality with parameter $\beta$ (PI($\beta$))  if there exists a
constant $C_{\mathrm{PI}}$ such that on each ball $B=B(z,r)$,
\begin{eqnarray*}
\int_B (f(x) - f_B)^2 d\mu(x) \le C_{\mathrm{PI}} r^{\beta} \int_B |\nabla f|^2 d\mu(x)
\end{eqnarray*}
 where $f_B=\frac{1}{\mu(B)} \int_B f(x)d\mu(x)$. 

We say that $(X,\mu,d)$ satisfies volume doubling (VD) if
there is a constant $C_{\mathrm{VD}}$ such that \begin{equation}
  V(x,2r)\le C_{\mathrm{VD}}V(x,r)\label{eq:VD}\end{equation}
for all $x\in X$ and $r>0.$  With
$\alpha=\log_{2}C_{\mathrm{VD}}$ and $r\ge s,$ we can use the equivalent formulation to compare volumes of nearby balls \begin{equation}
  \frac{V(x,r)}{V(y,s)}\le
  C_{\mathrm{VD}}\left(\frac{d(x,y)+r}{s}\right)^{\alpha}.\label{eq:vd1}\end{equation}

\section{Results and Discussion}\label{Results}
For the setting presented, the main result is the following theorem.
\begin{thm} \label{main}
For any $\epsilon >0, t>0$, let $C_t$ denote the points in the sample path of the diffusion up to time $t$ and $C_t^{\epsilon}=\cup_{x\in C_t}B(x,\epsilon)$, the union of the radius $\epsilon$ balls centered at the points of  $C_t$.  Then for $\nu>0$ and the diffusion started from $x\in X$, 
\begin{equation}
\log \mathbb{E}^x\left[\exp\left(-\nu \mu\left(C_{t^{\beta}V(x,t)}^{\epsilon}\right)\right)\right] \simeq V(x,t).
\end{equation} 
\end{thm}

The lower bound portion of Theorem \ref{main} does not differ substantially from the classical approach \cite{DV}, \cite{Kac}, requiring only the careful estimation of the probability that the sausage remains inside of a ball of radius $r$ until at least time $r^{\beta}V(x,r)$.   To treat the more challenging upper bound result, we turn to a coarse graining method out of the field of random media \cite{SznBolt},\cite{Gibson}, \cite{SznBook}.  

To introduce this technique, for $\nu>0$ let $\mathbb{P}^{\nu}$ represent the law of the Poisson process in $X$ with rate function $\nu d\mu$ which is independent of the law $\mathbb{P}$ of the diffusion process.  The points $\{x_i\}$ of this process will represent the centers of hard obstacles $B(x_i,\epsilon)$ for some $\epsilon>0$.  In any realization of this Poisson process, the event that the diffusion path avoids the obstacles up to time $s$ is the same as the event that the Wiener sausage up to time $s$ contains no point of the Poisson process, the probability of which with respect to $\mathbb{P}^{\nu}$ is $\exp(-\nu \mu(C_s^{\epsilon}))$.  Averaging over all realizations of the point process, then, 
\begin{equation}
\mathbb{E}^x\left[\exp(-\nu \mu(C_s^\epsilon)) \right]= \mathbb{E}^{\nu}\left[\mathbb{P}^x\left[T>s\right]\right],
\end{equation}
where $T$ is the hitting time of the diffusion to the obstacle set $\{ B(x_i,\epsilon) \}_i$, and $\mathbb{E}^{\nu}$ is the expectation with respect to $\mathbb{P}^{\nu}$.  Determining the decay rate of this averaged survival probability is then accomplished by locating the balance between those obstacle configurations which are dense, in which the survival decays too quickly, and those  configurations which are themselves much less likely, but in which survival probability decays more slowly.  The proof appears in Section \ref{Proof}.

The results here show how the Gaussian estimates may be used to replace the parts of the proof which relied upon the structure of Euclidean space in the original demonstration of this technique in \cite{SznBolt}.  This type of use of the Gaussian estimates has also been used in \cite{LSCBendikov} to study the properties of sample paths of diffusions in local Dirichlet spaces.  In future work, the authors hope to determine whether or not these techniques might allow for the more involved coarse graining technique called the method of enlargement \cite{SznBook} to be used in the local Dirichlet space context to obtain more precise asymptotic information.  

\subsection{Examples}\label{Examples}
Although the case of the translation invariant spaces $\Z^d$ \cite{DV}, $\R^d$ \cite{SznSausage}, and Cayley graphs of finitely generated polynomial growth groups \cite{Erschler2}, \cite{Erschler}, as well as the potentially translation invariant infinite graphs with Gaussian and sub-Gaussian estimates \cite{Gibson} have been studied, the setting of a local Dirichlet space broadens the class of examples that one can consider to a large class of continuous spaces without explicit symmetries.
\subsubsection{$\beta=2$}\label{Beta2}
A particularly nice class of examples of a local Dirichlet space with Gaussian bounds are complete Riemannian manifolds with non-negative Ricci curvature that satisfy the volume estimates (Eqs. (\ref{gb1}) and (\ref{eq:volcondbeta})).  The manifold structure is an example of a local Dirichlet space, and the non-negative Ricci curvature guarantees the Gaussian bounds.  Note that any Ahlfors regular manifold ($V(x,r)\simeq r^{\gamma}$ for some $\gamma>0$) will satisfy the volume estimates, although in general Ahlfors regularity is not a necessary condition. 

A Riemannian complex with an underlying finitely generated polynomial growth group will satisfy Gaussian bounds on the heat kernel.  (See \cite{PivarskiLSC} for a detailed description.)  For any polynomial volume growth finitely generated group, $G$, it is possible to construct many such objects.   A simple example is formed by taking the Cartesian product of a Cayley graph of $G$ with  $ [0,1]$ where $(g,0)$ and $(h,1)$ are identified whenever $hg^{-1}$ is an element of the generating set.  

\subsubsection{$\beta <2$}\label{Betalessthan2}
Fractals form a natural class of examples where $\beta$-Gaussian bounds hold \cite{Barlow}.  In \cite{P-P}, Pietruska-Pa{\l}uba proved more precise Wiener sausage asymptotics for simple nested fractals.  The local Dirichlet space setting allows us to extend our result to a broader class of fractal-like objects.

In \cite{BBK}, Barlow et al. show that if $X$ and $Y$ are metric measure spaces (connected locally compact complete separable metric space with geodesic measure) that have regular, strong local Dirichlet forms on $L^2(X, \mu)$, are stochastically complete, and that are roughly isometric to one another,  then $X$ satisfies GE($\beta$) if and only if $Y$ also satisfies GE($\beta$).
This allows one to construct a number of examples which are are roughly isometric to fractals with known bounds.  One such example is a tube-like fractal manifold where the edges of a fractal graph are replaced by tubes which are identified at the vertices \cite{BBK}; another is a cable system of a graph formed by replacing the edges of a graph with the interval $[0,1]$ where vertices are appropriately identified.  

For certain fractal manifolds, such as a pre-fractal generalized Sierpinski carpet, the $\beta$-Gaussian estimates hold only for $t\ge d(x,y)$ \cite{BarlowBassBrownian}.  In this situation, the proof of  Lemma \ref{lem:LBlemma} does not apply.  It would be interesting to see if the Wiener sausage asymptotics hold for spaces with these estimates.

\section{Proof of Theorem}\label{Proof}
\subsection{Lower Bound}\label{Lower}
Let $\tau_{B(x,r)}$ be the exit time of the diffusion process starting from $x$ from the ball $B(x,r)$.   Let $s=t^{\beta}V(x,t)$.  Then with $t>2\epsilon$,
\begin{eqnarray}
\mathbb{E}^x[\exp(-\nu  \mu(C_s^{\epsilon}))] & \ge & \mathbb{E}^x[\exp(-\nu V(x,t)); C_s^{\epsilon}\subset B(x,t)] \nonumber \\
 & \ge & \exp(-\nu V(x,t)) \mathbb{P}^x[\tau_{B(x,t-\epsilon)}>s]. \label{RHSLB}
\end{eqnarray}
Lemma \ref{lem:LBlemma} below demonstrates that for some constants $c, C', A>0$ 
\begin{equation*}
\mathbb{P}^x[\tau_{B(x,t-\epsilon)}>s]\ge \frac{c e^{-C'(t-\epsilon)}}{2V(x,(t-\epsilon)^{1/\beta})}\exp(-s\lambda(B(x,At-A\epsilon))). 
\end{equation*}
By Lemma 5.12 of \cite{LSCHebisch}, $\lambda(B(x,At-A\epsilon))\ge c/t^{\beta}$, so that substituting $s=t^{\beta}V(x,t)$ we have 
$ \exp(-s\lambda(B(x,At-A\epsilon))) \ge \exp(-cV(x,t))$.
Finally, since $V(x,t)$ is at least linear in $t$, $\exp(-cV(x,t))$ is the dominant term in the right hand side of Eq. (\ref{RHSLB}), so that with a change of constants and sufficiently large $t$, 
\begin{equation*}
\mathbb{E}^x[\exp(-\nu  \mu(C_s^{\epsilon}))]  \ge  \exp((-\nu-c) V(x,t)),
\end{equation*}  which provides the lower bound of Theorem \ref{main}.  \qed 

\begin{lem}
\label{lem:LBlemma} There exist constants $c, A, C'$ such that for $\rho > \rho_0$  and for $\sigma$ sufficiently large with respect to $\rho,$ \[
\mathbb{P}^{x}\left[\tau_{B(x,\rho)}>\sigma\right]
\ge
\frac{c}{2V(x,\rho^{1/\beta})} e^{-C'\rho}
e^{-\sigma\lambda(B(x, A\rho))}.
\]

\end{lem}
\begin{proof} 
Let $B=B(x,\rho).$   We begin by considering a point $y$ which is in $B(x, A\rho)$ for some fixed $A<1$.  We apply first the Dynkin-Hunt formula for the killed heat kernel.  Then we bound the expectation above by the supremum of the kernel on the set.  Because $y \in B(x, A\rho)$, we have $d(y,x) \le A\rho$ and $d(y,\bdry B) \ge (1-A)\rho$.  We use the Gaussian upper and lower bounds on $h_t(x,y)$ to get an expression in terms of the volumes.  We then use volume doubling to remove the $y$ dependence. 
\begin{eqnarray*}
  h_{\rho}^B(x,y) 
&= &  
h_{\rho}(x,y) - \mathbb{E}^x[h_{\rho-\tau_B}(X_{\tau_B},y)1_{\tau_B\le \rho} ]
\\
&\ge &  
h_{\rho}(x,y) - \sup_{0<s<\rho}\sup_{z\in \del B} h_{s}(z,y)
\\
&\ge & 
\frac{c}{V(x,\rho^{1/\beta})}\exp\left(-C \rho A^{\frac{\beta}{\beta-1}}\right)
- 
\frac{C}{
V(y,\rho^{1/\beta})
}\exp\left(-c\rho(1-A)^{\frac{\beta}{\beta-1}}
\right)\\
&\ge & 
\frac{c}{V(x,\rho^{1/\beta})}\exp\left(- C \rho A^{\frac{\beta}{\beta-1}}\right)
- 
\frac{C C_{\mathrm{VD}}\left(\frac{A\rho+\rho^{1/\beta}}{\rho^{1/\beta}}\right)^{\alpha}}{
V(x,\rho^{1/\beta})
}\exp\left(-c\rho(1-A)^{\frac{\beta}{\beta-1}}
\right).
\end{eqnarray*}
We now have a lower bound on the heat kernel involving an expression depending on both the volume of a ball of radius $\rho^{1/\beta}$ and on an expression involving $\rho$ and $A$.  For a sufficiently large $\rho$ (to handle the constants in front of the exponentials) and a sufficiently small $A$ (to handle the constants in the exponentials), this expression will be positive. We want to bound:
\begin{eqnarray*}
c\exp\left(-C\rho A^{\frac{\beta}{\beta-1}}\right)
- 
C C_{\mathrm{VD}}\left(\frac{A\rho+\rho^{1/\beta}}{\rho^{1/\beta}}\right)^{\alpha}
\exp\left(-c\rho\left(1-A\right)^{\frac{\beta}{\beta-1}}
\right).
\end{eqnarray*}
We can make the exponentials comparable by setting $A= \frac{c^{(\beta-1)/\beta}}{(2C)^{(\beta-1)/\beta} +c^{(\beta-1)/\beta}}$.  
Rewriting these constants in terms of 
$C' = C A^{\beta/(\beta-1)} $ and a bit of algebra yields:
\begin{eqnarray*}
\exp\left(-C'\rho \right)
\left(
c
- 
C C_{\mathrm{VD}}
\left(
\left(\frac{C' \rho}{C}\right)^{\frac{\beta-1}{\beta}}
+1 \right)^{\alpha}
\exp\left(- C'\rho \right)
\right).
\end{eqnarray*}
Here, we can choose $\rho$ to be large enough for 
\begin{eqnarray*}
\frac{c}{2} \ge
C C_{\mathrm{VD}}
\left(
\left(\frac{C' \rho}{C}\right)^{\frac{\beta-1}{\beta}}
+1 \right)^{\alpha}
\exp\left(- C'\rho \right).
\end{eqnarray*}
For such a choice of $\rho$, we can combine these pieces to obtain for any  $y \in B(x. A\rho)$:
\begin{eqnarray*}
  h_{\rho}^B(x,y) 
&\ge & 
\frac{c}{2V(x,\rho^{1/\beta})} \exp\left(-C'\rho \right).
\end{eqnarray*}
We will use that fact in the following computation:
\begin{eqnarray*}
  \mathbb{P}^{x}[\tau_{B}>\sigma] &=& \int_B h_{\sigma}^B(x,y)d\mu(y) \\
&\ge &  \int_{B(x, A\rho)} h_{\rho}^B(x,y)h_{\sigma-\rho}^B(y,y)d\mu(y)\\
&\ge  &  \int_{B(x, A\rho)} 
\frac{c}{2V(x,\rho^{1/\beta})} e^{-C'\rho}
h_{\sigma-\rho}^{B(x, A\rho)}(y,y)d\mu(y)\\
&=  & 
\frac{c}{2V(x,\rho^{1/\beta})} e^{-C'\rho}
\Tr(H_{\sigma-\rho}^{B(x, A\rho)})\\
&\ge & 
\frac{c}{2V(x,\rho^{1/\beta})} e^{-C'\rho}
e^{-(\sigma - \rho)\lambda(B(x, A\rho))}\\
&\ge & 
\frac{c}{2V(x,\rho^{1/\beta})} e^{-C'\rho}
e^{-\sigma\lambda(B(x, A\rho))}.
\end{eqnarray*}
The first inequalities make use of semigroup properties and GE($\beta$).
The next-to-last inequality comes from the fact that since
$H^{B(x, A\rho)}_{\sigma-\rho}$ has non-negative eigenvalues based on the
eigenvalues of $\Delta^{B(x, A\rho)}$, we can bound the trace below by
$e^{-(\sigma - \rho)\lambda(B(x, A\rho))}$. 
\qed
\end{proof}

\subsection{Upper Bound}\label{Upper}
\begin{proof}
As described in Section \ref{Results}, let $\mathbb{P}^{\nu}$ represent the law of the Poisson process in $X$ with rate function $\nu d\mu$ taken to be independent of the law $\mathbb{P}^x$ of the diffusion process starting at $x$.  
Let $\Omega$ be the probability space of the Poisson point process.  By $\{x^{\omega}_i\}$ we mean the points of the process for particular realization $\omega\in\Omega$.  These points will represent the centers of hard obstacles $B(x^{\omega}_i,\epsilon)$, and $T$ is the first time that the diffusion process enters the collection of obstacles. At times it is convenient to refer to $\cup_i B(x^{\omega}_i,\epsilon)$ by the symbol $\U$. 
Continuing with $s=t^{\beta}V(x,t)$, for any $N>0$, for a particular realization of the Poisson process,  
\begin{equation}\label{eq:inorout}
\mathbb{P}^{x}\left[T>s\right] \le \mathbb{P}^{x}\left[T\wedge \tau_{B\left(x,Ns\right)}>s\right]+\mathbb{P}^x\left[\tau_{B\left(x,Ns\right)}<s\right].
\end{equation}
By Lemma \ref{lem:GaussianExercise}, 
\begin{eqnarray*}
 \mathbb{P}^x\left[\tau_{B\left(x,Ns\right)}<s\right] &\le & C\exp\left(-c\left(\frac{(Ns)^{\beta}}{s}\right)^{1/(\beta-1)}\right) \\
  & \le & C\exp\left(-cN^{\beta/\beta-1}\left(s^{\beta-1}\right)^{1/(\beta-1)}\right) \\
   & \le & C\exp(-cs)=C\exp(-ct^{\beta}V(x,t)).
 \end{eqnarray*} 
 As demonstrated in Lemma (\ref{lem:GaussianExercise}), this decays faster than $\mathbb{P}^{x}\left[T\wedge \tau_{B\left(x,Ns\right)}>s\right]$, and will not contribute to the log-asymptotic behavior.

For convenience, set $\tilde{T}=T\wedge \tau_{B\left(x,Ns\right)}$, and let $B_{r}^{\omega}=B(x,r)\setminus \U$. Now using H\"{o}lder's inequality and the spectral theorem for $H^{B_{Ns}^{\omega}}_{s}$, 
\begin{eqnarray*}
\mathbb{P}^{x}\left[\tilde{T}>s\right] 
 & = & \left\langle \delta_{x},H^{B_{Ns}^{\omega}}_{s}\mathbf{1}_{B_{Ns}^{\omega}}\right\rangle \\
 & \le & \exp\left(-s\lambda\left(B_{Ns}^{\omega}\right)\right)\left\Vert \delta_{x}\right\Vert _{2}\left\Vert \mathbf{1}_{B_{Ns}^{\omega}}\right\Vert _{2}.
 \end{eqnarray*}
Averaging over $\Omega$, 
 \begin{equation}
 \mathbb{E}^{\nu}\left[\mathbb{P}^x\left[\tilde{T}>s\right]\right] \le \sqrt{V(x,Ns)}\mathbb{E}^{\nu}\left[\exp\left(-s\lambda\left(B_{Ns}^{\omega}\right)\right)\right].
 \end{equation}
Since $V(x,t)$  is subexponential (due to volume doubling), it suffices to show that
\begin{equation}
\limsup_{t\rightarrow\infty}V(x,t)^{-1}\log\mathbb{E}^{\nu}\left[\exp\left(-s\lambda\left(B_{Ns}^{\omega}\right)\right)\right]  <  0.
\end{equation} 
Recalling the definition of $s$, this automatically holds whenever $\lambda\left(B_{Ns}^{\omega}\right)\ge c_0/t^{\beta}$ for any constant $c_0>0$.  Therefore we restrict to the complementary case, reserving the right to choose $c_0$ at a later point.  Thus, as
\begin{equation}
\mathbb{E}^{\nu}\left[\exp\left(-s\lambda\left(B_{Ns}^{\omega}\right)\right) \mathbf{1}_{\left\{ \lambda\left(B_{Ns}^{\omega}\right) \le c_{0}t^{-\beta}\right\} }\right]  \le  \mathbb{P}^{\nu}\left[\lambda\left(B_{Ns}^{\omega}\right)\leq c_{0}t^{-\beta}\right],
\end{equation}
it now suffices to show that 
\begin{equation}
\limsup_{t\rightarrow\infty}V(x,t)^{-1}\log\mathbb{P}^{\nu}\left[\lambda\left(B_{Ns}^{\omega}\right)\leq c_{0}t^{-\beta}\right]  <  0.
\end{equation}
Fix $t>0$ and a realization $\omega$ of the Poisson point process.
Define an $t$-net $\{K_i,k_i \}$ on $X$.  This a family of balls $K_i=B(k_i,t)$ for which $X\subset \cup_i K_i$, but $\{B(k_i,t/2)\}$ form a pairwise disjoint collection of sets.  Since the volume growth is at least linear, by volume doubling the number of net elements which overlap with $B(x,s)$ is at most a polynomial $Q(s)$.  Also, every point of $y \in X$ may be contained in at most 
\begin{equation}
C_{\mathrm{over}}=C_{VD}\left(\frac{t+2t}{t/2}\right)^{\alpha}\ge \frac{V(y,2t)}{V(z,t/2)}
\end{equation}
$t$-net elements, where $z\in B(y,t)$, by Eq. (\ref{eq:vd1}).
We can now relate $\lambda\left(B_{Ns}^{\omega}\right)$ to the proportion of space taken up by the obstacles inside the $t$-net elements.
Consider $f\in L^2(B \setminus \U), \supp(f) \subset B \setminus \U, ||f||_2 \ne 0, f\in\Dom(\E)$.  Although not explicitly stated, if $f$ is identically 0 on a $K_i$, then that $K_i$ is left out of the intermediate steps.
 Now using Lemma \ref{Thirring's inequality} and PI($\beta$),
\begin{align*}
\lefteqn{\int_X |\grad f(x)|^2 d\mu(x)}\\
 & \ge C_{\mathrm{over}}^{-1}\sum_{i=1}^{Q(t)}\left(
\int_{K_i} |\grad f(x)|^2 d\mu(x)
\right)\\
& \ge C_{\mathrm{over}}^{-1}\sum_{i=1}^{Q(t)}\left( 
\int_{K_i} |f(x) |^2 d\mu(x) 
\inf_{\stackrel{f\in L^2(K_i \setminus \U)}{\stackrel{ \supp(f) \subset K_i \setminus \U}{||f||_2 \ne 0}}}
\frac{\int_{K_i} |\grad f(x)|^2 d\mu(x)}{\int_{K_i} |f(x)|^2 d\mu(x)}
\right)
\\
& \ge C_{\mathrm{over}}^{-1}\sum_{i=1}^{Q(t)}\left( 
\int_{K_i} |f(x)|^2 d\mu(x) 
\inf_{\stackrel{f\in L^2(K_i)}{\stackrel{ \supp(f) \subset K_i}{ ||f||_2 \ne 0}}} 
\frac{ \int_{K_i} |\grad f(x)|^2 d\mu(x)}{\int_{K_i} |f(x)|^2 d\mu(x)} 
\frac{\mu\left(\U\cap K_{i} \right)}{\mu\left(K_{i}\right)}
\right)
\\
& \ge C_{\mathrm{over}}^{-1}\sum_{i=1}^{Q(t)}\left( 
\int_{K_i} |f(x)|^2 d\mu(x) 
\inf_{\stackrel{f\in L^2(K_i)}{\stackrel{ \supp(f) \subset K_i}{ ||f||_2 \ne 0}}}
\frac{C_{PI}t^{-\beta} \int_{K_i} |f(x)|^2 d\mu(x)}{\int_{K_i} |f(x)|^2 d\mu(x)} 
\frac{\mu\left(\U\cap K_{i}\right)}{\mu\left(K_{i}\right)}
\right)
\\
& \ge C_{\mathrm{over}}^{-1} C_{PI}t^{-\beta} 
\inf_{i \in [1, Q(t)]}
\frac{\mu\left( \U\cap K_{i} \right)}{\mu\left(K_{i}\right)}
\int_{B} |f(x)|^2 d\mu(x) 
.
\end{align*}
Therefore, (absorbing the constants)
\begin{align*}
\lambda\left(B_{Ns}^{\omega}\right) & \ge ct^{-\beta}\inf_{i\in[1,Q(s)]}\left\{ \frac{\mu\left( \U\cap K_{i} \right)}{\mu\left(K_{i}\right)}\right\} \end{align*}
and 
\begin{align}
\mathbb{P}^{\nu}\left[\lambda\left(B_{Ns}^{\omega}\right)\leq c_{0}t^{-\beta}\right] & \le\mathbb{P}^{\nu}\left[\inf_{i\in[1,Q(s)]}\frac{\mu\left(\U\cap K_{i} \right)}{\mu\left(K_{i}\right)}\leq c^{-1}c_{0}\right]\\
 & \le Q(s)\sup_{i\in[1,Q(s)]}\mathbb{P}^{\nu}\left[\frac{\mu\left( \U\cap K_{i} \right)}{\mu\left(K_{i}\right)}\leq c^{-1}c_{0}\right].\label{eq:withQ}\end{align}
Choosing $c_{0}$ sufficiently small with respect to $\nu$ 
and $c$, the right hand side of the above may be interpreted (after correcting the constant to take into account the possibility of multiple Poisson arrivals in an $\epsilon$-ball in $K_i$) as the event that there are proportionally less Poisson arrivals in $K_i$ than the expected number $\nu\mu(K_i)$.  Then, as the volume of the radius $t$ balls grows with $t$, it follows from a Cramer type estimate (for example, see Theorem 2.2.1 of \cite{Dembo}) that for some constant $c'$ and all large enough $t$,
\begin{eqnarray*}
\mathbb{P}^{\nu}\left[\frac{\mu\left(\U\cap B\left(k_{i},t/2\right)\right)}{\mu\left(B\left(k_{i},t/2\right)\right)}\leq c^{-1}c_{0}\right] & \leq & \exp\left(-c'\mu\left(K_{i}\right)\right).
\end{eqnarray*}
Since $Q(s)$ is polynomial in $s$ (and in $t$ as well), and since $V(x,t)$ is at least linear in $t$, the exponential term in the right hand side of Eq. (\ref{eq:withQ}) dominates as $t$ grows large. This implies  
\begin{align*}
\lefteqn{\limsup_{t\to\infty}V(x,t)^{-1}\log\mathbb{P}^{\nu}\left[\lambda\left(B_{s}^{\omega}\right)\leq c_{0}t^{-\beta}\right]}\\
 & <-c'\liminf_{t\to\infty}V(x,t)^{-1}\inf_{i\in[1,Q(s)]}\mu\left(K_{i}\right).
 \end{align*}
 Whenever the right had side of the above is strictly less than zero, the upper bound of Theorem \ref{main} holds.  But this is equivalent to the volume regularity requirement Eq. (\ref{eq:volcondbeta}).
\end{proof}

 We conclude by establishing the two lemmas used above.  In the lemma below, the Gaussian estimates are used to demonstrate that 
 $$\mathbb{P}^x\left[\tau_{B\left(x,Ns\right)}<s\right] \le \exp(-cs).$$ A result of this type appears in \cite{LscHebisch0}, and an exact discrete analogue to the below in \cite{Gibson2}.
 \begin{lem}\label{lem:GaussianExercise}
   Suppose that \emph{GE($\beta$)}
holds. Then for the diffusion process starting from $y\in X$, \[
\mathbb{P}^{y}\left[\sup_{1\le s\le t}d\left(y,X_{s}\right)\ge\sigma\right]\le C\exp\left(-c\left(\frac{\sigma^{\beta}}{t}\right)^{1/(\beta-1)}\right).\]
\end{lem}
\begin{proof}
Note that the result is void if $\sigma^{\beta}/t$ is small. Let
$B^{j}=B\left(y,2^{j}\sigma\right)$ be the ball of radius $2^{j}\sigma$
about $y$. We first prove that \begin{equation}
\mathbb{P}^{y}\left[X_{t}\notin B^{0}\right]\leq C\exp\left(-c\left(\frac{\sigma^{\beta}}{t}\right)^{1/(\beta-1)}\right).\label{eq:escapedecay}\end{equation}
 Indeed, $\mathbb{P}^{y}\left[X_{t}\notin B^{0}\right]=\int_{X\setminus B^{0}}h_{t}(y,x)d\mu(x)$,
so that by Eqs. (\ref{eq:GUE}), (\ref{eq:vd1}), and a change of constants
we obtain\begin{align*}
\lefteqn{\int_{x\not\in B(y,\sigma)}h_{t}(y,x)d\mu(x)}\\
\leq & \frac{C}{V\left(y,t^{1/\beta}\right)}\int_{X\setminus B^{0}}\exp\left(-c\left(\frac{(d(y,x))^{\beta}}{t}\right)^{\frac{1}{\beta-1}}\right)d\mu(x)\\
= & \frac{C}{V\left(y,t^{1/\beta}\right)}\sum_{j=1}^{\infty}\left(\int_{B^{j}\setminus B^{j-1}}\exp\left(-c\left(\frac{(d(y,x))^{\beta}}{t}\right)^{\frac{1}{\beta-1}}\right)\right)d\mu(x)\\
\leq & C\sum_{j=1}^{\infty}\frac{V\left(y,2^{j}\sigma\right)}{V\left(y,t^{1/\beta}\right)}\exp\left(-c\left(\frac{\left(2^{(j-1)}\sigma\right)^{\beta}}{t}\right)^{\frac{1}{\beta-1}}\right)\\
\leq & C\left(\frac{\sigma}{t^{1/\beta}}\right)^{\alpha}\exp\left(-c\left(\frac{\sigma^{\beta}}{2t}\right)^{\frac{1}{\beta-1}}\right)\sum_{j=1}^{\infty}C_{\mathrm{VD}}^{j}\exp\left(-c\left(\frac{\sigma^{\beta}}{2t}\right)^{\frac{1}{\beta-1}}\left(2^{\frac{\beta}{\beta-1}\left(j\right)}-1\right)\right)\\
\leq & C\exp\left(-c\left(\frac{\sigma^{\beta}}{t}\right)^{\frac{1}{\beta-1}}\right).\end{align*}
The necessary size of $\sigma^{\beta}/t$ only depends on the size
of $c,C$ from the Gaussian-type estimates and $C_{\mathrm{VD}}$ from volume doubling.

Let $L_{\sigma}=\inf\left\{ t\left|d\left(X_{0},X_{t}\right)>\sigma\right.\right\} .$
Intersecting the event $\left\{ X_{t}\notin B(y,\sigma/2)\right\} $
by the event $\left\{ L_{\sigma}\leq t\right\} $ and using the strong
Markov property, \begin{eqnarray*}
\lefteqn{\mathbb{P}^{y}\left[X_{t}\notin B(y,\sigma/2)\right]}\\
 & \geq & \mathbb{P}^{y}\left[X_{t}\notin B(y,\sigma/2),L_{\sigma}\leq t\right]\\
 & \geq & \mathbb{P}^{y}\left[L_{\sigma}\leq t\right]-\mathbb{P}^{y}\left[X_{t}\in B(y,\sigma/2),L_{\sigma}\leq t\right]\\
 & \geq & \mathbb{P}^{y}\left[L_{\sigma}\leq t\right]-\mathbb{E}^{y}\left[\mathbb{P}^{X_{L_{\sigma}}}\left[X_{t-L_{\sigma}}\notin B\left(X_{L_{\sigma}},\sigma/2\right)\right]\mathbf{1}_{\left\{ L_{\sigma}\leq t\right\} }\right].\end{eqnarray*}
 When $\sigma^{\beta}/t$ is large enough, Eq. (\ref{eq:escapedecay})
may be applied to obtain \[
\sup_{x,k\leq t}\mathbb{P}^{x}\left[X_{k}\notin B(x,\sigma/2)\right]\leq C'<1.\]
 Combining this with Eq. (\ref{eq:escapedecay}) once more shows that
\begin{eqnarray*}
\mathbb{P}^{y}\left[L_{\sigma}\leq t\right](1-C') & \leq & \mathbb{P}^{y}\left[X_{t}\notin B(y,\sigma/2)\right]\\
 & \leq & C\exp\left(-c\left(\frac{\sigma^{\beta}}{t}\right)^{1/(\beta-1)}\right).\end{eqnarray*}

 \end{proof}
  
 The following lemma (nearly identical to Lemma 3.3 of \cite{SznBolt})
is used in the upper bound proof to transfer the focus from the eigenvalue of $B_{NS}^{\omega}$ to the proportion of volume taken up by obstacles in the $t$-net balls.

\begin{lem}
\label{Thirring's inequality} For a compact subset $U\subset X$
and $A\subset U$, define \[
\lambda_{A}(U)=\inf_{\stackrel{{\scriptstyle f\in L^{2}(U,\mu)}}{\|f\|_{2}^{2}=1}}\left\{ \mathcal{E}(f,f)|\mathrm{supp}f\subset U\setminus A\right\} .\]
 Then \[
\lambda_{A}(U)\ge\lambda(U)\frac{\mu(A\cap U)}{\mu(U)}.\]
 
\end{lem}
\begin{proof}
Choose $f\in L^{2}\left(U,\mu \right) \cap \Dom(\E)$ to be vanishing on $A.$ Let
$g=f-\langle f,\psi\rangle \psi$ where $\psi=\mathbf{1}_{U}/\sqrt{\mu(U)}$.  Note that both $\langle g,\psi\rangle$ and $\mathcal{E}(\psi,\psi)$ are 0.  With $f$ vanishing on $A,$ for any $\alpha,\beta>0$, 
\begin{eqnarray*}
 \mathcal{E}(f,f) + \alpha||f||_2^2 
&=& \mathcal{E}(g+\langle f,\psi\rangle \psi,g+\langle f,\psi\rangle \psi) + ||(\alpha +\mathbf{1}_{U}\beta) f||_2^2 \\
&=& \mathcal{E}(g,g) + ||(\alpha +\mathbf{1}_{U}\beta) f||_2^2 .
 \end{eqnarray*}
 By H\"{o}lder's inequality, \begin{eqnarray}
\langle f,\psi\rangle^{2}
 & \le & ||\left(\alpha+\beta\mathbf{1}_{A}(\cdot)\right)^{- \frac{1}{2}}\psi||_2^2
  ||\left(\alpha+\beta\mathbf{1}_{A}(\cdot)\right)^{\frac{1}{2}} f||_2^{2} \label{eq:cauchy} 
  \end{eqnarray}
 Since $\langle g,\psi\rangle=0$ , \[\mathcal{E}(g,g) \ge\lambda(U)\left\Vert g\right\Vert _{2}^{2},\]
 which with Eq. (\ref{eq:cauchy}) gives
 \begin{eqnarray*}
 \mathcal{E}(g,g) +||\left(\alpha+\beta\mathbf{1}_{A}(\cdot)\right)^{\frac{1}{2}} f||_2^{2}
 & \geq & 
 \lambda(U)\left\Vert g\right\Vert _{2}^{2} + \left(  ||\left(\alpha+\beta\mathbf{1}_{A}(\cdot)\right)^{- \frac{1}{2}}\psi||_2^2
 \right)^{-1} \langle f,\psi\rangle^{2} .
 \end{eqnarray*}
 Noticing that both $||g||_2^2$ and $\langle f,\psi\rangle^{2}$ are less than $||f||_2^2$
it follows that 
\begin{eqnarray}
\frac{ \mathcal{E}(f,f) + \alpha||f||_2^2}{||f||_2^2} 
&\ge& \lambda(U)\frac{|| g||_{2}^{2}}{||f||_2^2} + \left(  ||\left(\alpha+\beta\mathbf{1}_{A}(\cdot) \right)^{- \frac{1}{2}}\psi||_2^2
 \right)^{-1} \frac{\langle f,\psi\rangle^{2}}{||f||_2^2} \nonumber \\
 & \ge & \min\left(\lambda(U), \left(  ||\left( \alpha+\beta\mathbf{1}_{A}(\cdot) \right)^{- \frac{1}{2}}\psi||_2^2 \right)^{-1}  \right) \label{eq:lincomb}
 \end{eqnarray}
 Because $\psi = \mathbf{1}_{U}/\sqrt{\mu(U)}$, as $\beta$ goes to infinity,\[
||\left(\frac{1}{\alpha+\beta\mathbf{1}_{A}(\cdot)}\right)^{\frac{1}{2}}\psi||_2^2
\rightarrow\frac{\mu\left(U\setminus A\right)}{\alpha \mu(U)},\]
 which with Eq. (\ref{eq:lincomb}) shows that\[
\lambda_{A}(U)\geq\min\left\{ \mu(U),\alpha\frac{\mu(U)}{\mu\left(U\setminus A\right)}\right\} -\alpha.\]
 Taking $\alpha=\lambda(U)\frac{\mu\left(U\setminus A\right)}{\mu(U)}$, completes
the proof.
\end{proof}

The authors would like to thank Laurent Saloff-Coste for helpful discussions and inspiration.

\end{document}